\documentclass{amsart}
\usepackage{graphicx}
\usepackage{latexsym,amsmath,amssymb}

\theoremstyle{definition}
\theoremstyle{remark}

\numberwithin{equation}{section}

\begin{document}

\title[Weighted conditional type operators ]
{ A classes of operators with normal Aluthge transformation }

\author{\sc\bf Y. Estaremi  }
\address{\sc Y. Estaremi  }
\email{yestaremi@pnu.ac.ir}

\address{Department of mathematics, Payame Noor university , p. o. box: 19395-3697, Tehran,
Iran}

\thanks{}

\thanks{}



\date{}

\dedicatory{}

\commby{}

\begin{abstract}
In this paper, we show that the generalized Aluthge
transformations of a large class of operators (weighted
conditional type operators) are normal. As a consequence, the
operator $M_wEM_u$ is $p$-hyponormal if and only if it is normal,
and under a weak condition, if $M_wEM_u$ is normal, then the
Holder inequality turn into equality for $w,u$. Also, we conclude
that, invertible weighted conditional type operators are normal.
In the end we give some applications of $p$-hyponormal weighted
conditional type operators. In the end, some examples are provided
to illustrate concrete application of the main results of the
paper.

\noindent {}
\end{abstract}

\maketitle

\section{ \sc\bf Introduction }
As is well-known, conditional expectation operators on various
function space exhibit a number of remarkable properties related
either to the underlying order structure, or to the metric
structure when the function space is equipped with a norm. Such
operators are necessarily positive projections which are averaging
in a precise sense to be described below and in certain normed
function space are contractive for the given norm. Theory of
multiplication conditional type operators is one of important
arguments in the connection of operator theory and measure theory.
Conditional expectations have been studied in an operator
theoretic setting, by, for example, De pagter and Grobler \cite{g}
and Rao \cite{rao1, rao2}, as positive operators acting on
$L^p$-spaces or Banach function spaces. In \cite{mo}, S.-T. C. Moy
characterized all operators on $L^p$ of the form $f\rightarrow
E(fg)$ for $g$ in $L^q$ with $E(|g|)$ bounded.  Also, some results
about multiplication conditional expectation operators can be
found in \cite{dou, her, lam}. In \cite{dhd} P.G. Dodds, C.B.
Huijsmans and B. de Pagter showed that lots of operators are of
the form of weighted conditional type operators. The class of
weighted conditional type operators contains composition
operators, multiplication operators, weighted composition
operators, some integral type operators, ... . These are some
reasons that stimulate us to consider weighted conditional type
operators in our work. So, in \cite{e,ej, ej2} we investigated
some classic properties of multiplication conditional expectation
operators $M_wEM_u$ on $L^p$ spaces. In this paper we state that
the
Aluthge transformations of weighted conditional type operators are normal, so we use this property and give some applications.\\

\section{ \sc\bf Preliminaries}
Let $(X,\Sigma,\mu)$ be a complete $\sigma$-finite measure space.
For any sub-$\sigma$-finite algebra $\mathcal{A}\subseteq
 \Sigma$, the $L^2$-space
$L^2(X,\mathcal{A},\mu_{\mid_{\mathcal{A}}})$ is abbreviated  by
$L^2(\mathcal{A})$, and its norm is denoted by $\|.\|_2$. All
comparisons between two functions or two sets are to be
interpreted as holding up to a $\mu$-null set. The support of a
measurable function $f$ is defined as $S(f)=\{x\in X; f(x)\neq
0\}$. We denote the vector space of all equivalence classes of
almost everywhere finite valued measurable functions on $X$ by
$L^0(\Sigma)$.

\vspace*{0.3cm} For a sub-$\sigma$-finite algebra
$\mathcal{A}\subseteq\Sigma$, the conditional expectation operator
associated with $\mathcal{A}$ is the mapping $f\rightarrow
E^{\mathcal{A}}f$, defined for all non-negative measurable
function $f$ as well as for all $f\in L^2(\Sigma)$, where
$E^{\mathcal{A}}f$, by the Radon-Nikodym theorem, is the unique
$\mathcal{A}$-measurable function satisfying
$\int_{A}fd\mu=\int_{A}E^{\mathcal{A}}fd\mu, \ \ \ \forall A\in
\mathcal{A}$. As an operator on $L^{2}({\Sigma})$,
$E^{\mathcal{A}}$ is idempotent and
$E^{\mathcal{A}}(L^2(\Sigma))=L^2(\mathcal{A})$. This operator
will play a major role in our work. Let $f\in L^0(\Sigma)$, then
$f$ is said to be conditionable with respect to $E$ if
$f\in\mathcal{D}(E):=\{g\in L^0(\Sigma): E(|g|)\in
L^0(\mathcal{A})\}$. Throughout this paper we take $u$ and $w$ in
$\mathcal{D}(E)$. If there is no possibility of confusion, we
write $E(f)$ in place of $E^{\mathcal{A}}(f)$.  A detailed
discussion about this operator may be found in \cite{rao}. \\

Let $\mathcal{H}$ be the infinite dimensional complex Hilbert
space and let $\mathcal{L(H)}$ be the algebra of all bounded
operators on $\mathcal{H}$. It is known that an operator $T$ on a
Hilbert space is $p$-hyponormal if $(T^{\ast}T)^p\geq
(TT^{\ast})^p$, for $0<p<\infty$. Every operator $T$ on a Hilbert
space $\mathcal{H}$ can be decomposed into $T = U|T|$ with a
partial isometry $U$, where $|T| = (T^*T)^{\frac{1}{2}}$ . $U$ is
determined uniquely by the kernel condition $\mathcal{N}(U) =
\mathcal{N}(|T|)$. Then this decomposition is called the polar
decomposition. The Aluthge transformation $\widehat{T}$ of the
operator $T$ is defined by
$\widehat{T}=|T|^{\frac{1}{2}}U|T|^{\frac{1}{2}}$. The operator
$T$ is said to be positive operator and write $T\geq 0$, if
$\langle
Th, h\rangle\geq 0$, for all $h\in \mathcal{H}$. \\
In this paper we will be concerned with normality,
$p$-hyponormality and invertibility of weighted conditional type
operators.

\section{ \sc\bf Main results }

In the first we reminisce some properties of weighted conditional
type operators, that we have proved in
\cite{ej}.\\

 The operator $T=M_wEM_u$ is bounded on $L^{2}(\Sigma)$ if and
only if $(E|w|^{2})^{\frac{1}{2}}(E|u|^{2})^{\frac{1}{2}} \in
L^{\infty}(\mathcal{A})$, and in this case its norm is given by
$\|T\|=\|(E(|w|^{2}))^{\frac{1}{2}}(E(|u|^{2}))^{\frac{1}{2}}\|_{\infty}$.
The unique polar decomposition of bounded operator $T=M_wEM_u$ is
$U|T|$, where

$$|T|(f)=\left(\frac{E(|w|^{2})}{E(|u|^{2})}\right)^{\frac{1}{2}}\chi_{S}\bar{u}E(uf)$$
and
 $$U(f)=\left(\frac{\chi_{S\cap
 G}}{E(|w|^{2})E(|u|^{2})}\right)^{\frac{1}{2}}wE(uf),$$
for all $f\in L^{2}(\Sigma)$, where $S=S(E(|u|^2))$ and
$G=S(E(|w|^2))$. Also, the Aluthge transformation of $T=M_wEM_u$
is
$$\widehat{T}(f)=\frac{\chi_{S}E(uw)}{E(|u|^{2})}\bar{u}E(uf), \ \ \ \ \ \  \ \ \ \  \ \  \ f\in L^{2}(\Sigma).$$
For each $\epsilon>0$ we have
$$\widehat{T}_{\epsilon}(f)=|T|^{\epsilon}U|T|^{1-\epsilon}(f)=(E(|u|^2))^{-1}\chi_{S}E(uw)\overline{u}E(uf)=\widehat{T},$$

$$\widehat{T}^{\ast}(f)=(E(|u|^2))^{-1}\chi_{S}\overline{uE(uw)}E(uf),$$ and
$$|\widehat{T}|(f)=|\widehat{T}^{\ast}|(f)=E(|u|^2))^{-1}\chi_{S}|E(uw)|\overline{u}E(uf).$$

By above information we have a nice conclusion about conditional
type operators $T=M_wEM_u$. For every conditional type operator
$T=M_wEM_u$, the Aluthge transformation $\widehat{T}$ is
normal.\\

Now we consider generalized Aluthge transformation of $T=M_wEM_u$.
Let $T=U|T|$ be the polar decomposition of $T=M_wEM_u$. For $r>0$
and $r\geq t\geq0$, let $\widetilde{T}=|T|^{t}U|T|^{r-t}$. By
[Lemma 3.3, \cite{ej}],  for all $f\in L^2(\Sigma)$ we have

$$\widetilde{T}(f)=E(|w|^2)^{\frac{r-1}{2}}E(|u|^2)^{\frac{r-3}{2}}\chi_{S\cap
G}E(uw)\bar{u}E(uf).$$ \\

\vspace*{0.3cm} {\bf Theorem 2.1.} For every multiplication
conditional type operator $T=M_wEM_u$, the operator
$\widetilde{T}$ is normal.

\vspace*{0.3cm} {\bf Proof }  Direct computations show that

$$\widetilde{T}^{\ast}\widetilde{T}(f)=E(|w|^2)^{r-1}E(|u|^2)^{r-2}\chi_{S\cap
G}|E(uw)|^2\bar{u}E(uf)=\widetilde{T}\widetilde{T}^{\ast}(f).$$

This means that $\widetilde{T}$ is normal.\\

\vspace*{0.3cm} {\bf Theorem 2.2.}  Let $T=M_wEM_u$ be
multiplication conditional type operator. Then $T$ is
$p$-hyponormal if and only if $T$ is normal.

\vspace*{0.3cm} {\bf Proof } By Theorem 2.1 and [Theorem 3, \cite{th}] we conclude that if $T$ is $p$-hyponormal, then $T$ is normal. The converse is clear.\\

\vspace*{0.3cm} {\bf Theorem 2.3.} If $T=M_wEM_u$ is invertible,
then it is normal.

\vspace*{0.3cm} {\bf Proof } The Theorem 3.1 of \cite{or}
guarantees that, if $T$ is invertible and $\widetilde{T}$ is
normal, then $T$ is normal. By Theorem 2.1 we conclude that if
$M_wEM_u$ is invertible, then it is normal.\\

\vspace*{0.3cm} {\bf Lemma 2.4.} Let $T=M_wEM_u$. Then
$|\widetilde{T}|=|\widetilde{T}^{\ast}|$.

\vspace*{0.3cm} {\bf Proof } Since $\widetilde{T}$ and
$\widetilde{T}^{\ast}$ are also weighted conditional type
operators, then we get that for all $f\in L^2(\Sigma)$

$$|\widetilde{T}|(f)=E(|w|^2)^{\frac{r-1}{2}}E(|u|^2)^{\frac{r-3}{2}}\chi_{S\cap
G}|E(uw)|\bar{u}E(uf)=|\widetilde{T}^{\ast}|(f).$$ \\

\vspace*{0.3cm} {\bf Theorem 2.5.} Let $T=M_wEM_u$ be
$p$-hyponormal. Then
$$|\widetilde{T}|=|\widetilde{T}^{\ast}|=|T|^r.$$

\vspace*{0.3cm} {\bf Proof } By using Lemma 2.4 and [Theorem 2,
\cite{th}] we have $|\widetilde{T}|=|\widetilde{T}^{\ast}|=|T|^r$.\\

\vspace*{0.3cm} {\bf Theorem 2.6.} If $T=M_wEM_u$ is
$p$-hyponormal, then $|E(uw)|^2=E(|u|^2)E(|w|^2)$ on $S(E(u))$.

\vspace*{0.3cm} {\bf Proof } Let $f\in L^2(\Sigma)$, by [Lemma
3.3, \cite{ej}] and Theorem 2.5, for all $f\in L^2(\Sigma)$ we
have

$$(E(|u|^2))^{\frac{r-2}{2}}(E(|w|^2))^{\frac{r}{2}}\chi_{S\cap G}\bar{u}E(uf)=(E(|u|^2))^{\frac{r-3}{2}}(E(|w|^2))^{\frac{r-1}{2}}\chi_{S\cap
G}|E(uw)|\bar{u}E(uf)$$ and so for $0<a\in L^2(\mathcal{A})$

$$((E(|u|^2))^{\frac{r-2}{2}}(E(|w|^2))^{\frac{r}{2}}-(E(|u|^2))^{\frac{r-3}{2}}(E(|w|^2))^{\frac{r-1}{2}}|E(uw)|)\chi_{S\cap
G}\bar{u}E(u)a=0.$$

By taking $E$ we have

$$((E(|u|^2))^{\frac{1}{2}}(E(|w|^2))^{\frac{1}{2}}-|E(uw)|)(E(|u|^2))^{\frac{r-3}{2}}(E(|w|^2))^{\frac{r-1}{2}}\chi_{S\cap
G}\bar{u}E(u)a=0.$$
And then
$$((E(|u|^2))^{\frac{1}{2}}(E(|w|^2))^{\frac{1}{2}}-|E(uw)|)|E(u)|^2=0.$$
This implies that $|E(uw)|^2=E(|u|^2)E(|w|^2)$ on $S(E(u))$.\\

\vspace*{0.3cm} {\bf Corollary 2.7.} If $T=M_wEM_u$ is
$p$-hyponormal and $S(E(u))=X$, then the conditional type Holder inequality for $u,w$ turns into equality, i.e, $|E(uw)|^2=E(|u|^2)E(|w|^2)$.\\

From now on, we shall denote by $\sigma_{p}(T)$, $\sigma_{jp}(T)$,
the point spectrum of $T$, the joint point spectrum of $T$,
respectively. A complex number $\lambda\in \mathbb{C}$ is said to
be in the point spectrum $\sigma_{p}(T)$ of the operator $T$, if
there is a unit vector $x$ satisfying $(T-\lambda)x=0$. If in
addition, $(T^{\ast}-\bar{\lambda})x=0$, then $\lambda$ is said to
be in the joint point spectrum $\sigma_{jp}(T)$ of $T$. If $A,
B\in \mathcal{B}(\mathcal{H})$, then it is well known that

$$\sigma_{p}(AB)\setminus\{0\}=\sigma_{p}(BA)\setminus\{0\}, \ \ \ \sigma_{jp}(AB)\setminus\{0\}=\sigma_{jp}(BA)\setminus\{0\}.$$\\

Let $A_{\lambda}=\{x\in X:E(u)(x)=\lambda\}$, for $0\neq\lambda\in
\mathbb{C}$. Suppose that $\mu(A_{\lambda})>0$. Since
$\mathcal{A}$ is $\sigma$-finite, there exists an
$\mathcal{A}$-measurable subset $B$ of $A_{\lambda}$ such that
$0<\mu(B)<\infty$, and $f=\chi_{B}\in L^p(\mathcal{A})\subseteq
L^p(\Sigma)$. Now
$$EM_u(f)-\lambda f=E(u)\chi_{B}-\lambda \chi_{B}=0.$$ This
implies that $\lambda\in \sigma_p(EM_u)$.\\
If there exists $f\in L^p(\Sigma)$ such that $f\chi_{C}\neq 0$
$\mu$-a.e, for $C\in \Sigma$ of positive measure and
$E(uf)=\lambda f$ for $0\neq \lambda \in \mathbb{C}$, then
$f=\frac{E(uf)}{\lambda}$, which means that $f$ is
$\mathcal{A}$-measurable. Therefore $E(uf)=E(u)f=\lambda f$ and
$(E(u)-\lambda)f=0$. This implies that $C\subseteq A_{\lambda}$
and so $\mu(A_{\lambda})>0$. Hence
$$\sigma_p(EM_u)=\{\lambda\in\mathbb{C}:\mu(A_{\lambda})>0\}.$$
Thus
$$\sigma_p(M_wEM_u)\setminus \{0\}=\{\lambda\in\mathbb{C}:\mu(A_{\lambda,w})>0\}\setminus \{0\},$$
where $A_{\lambda,w}=\{x\in X:E(uw)(x)=\lambda\}$.\\

\vspace*{0.3cm} {\bf Theorem 2.8.} If
$|E(uw)|^2=E(|u|^2)E(|w|^2)$, then
 $$\sigma_{p}(M_wEM_u)=\sigma_{jp}(M_wEM_u).$$

\vspace*{0.3cm} {\bf Proof.} Let $f\in L^2(\Sigma)\setminus \{0\}$
and $\lambda\in \mathbb{C}$, such that
 $wE(uf)=\lambda f$. Let $M=span\{f\}$, the closed linear subspace
 generated by $f$. Thus we can represent $T=M_wEM_u$ as the
 following $2\times 2$ operator matrix with respect to the
 decomposition $L^2(\Sigma)=M\oplus M^{\perp}$,
 $$
 T= \left[
         \begin{array}{rr}
              M_{\lambda} & PM_{w}EM_{u}- M_{\lambda} \\
              0 &  M_{w}EM_{u}-PM_{w}EM_{u}
          \end{array} \right].
$$
where $P$ is is the orthogonal projection of $L^2(\Sigma)$ onto
$M$. Since $|E(uw)|^2=E(|u|^2)E(|w|^2)$, then we have
$P(|T^2|-|T^{\ast}|^2)P\geq0$. Direct computation shows that
$$P|T^2|^2P= \left[
         \begin{array}{rr}
              M_{|\lambda|^4} & 0 \\
              0 &  0
          \end{array} \right]
.$$

Thus
$$\left[
         \begin{array}{rr}
              M_{|\lambda|^2} & 0 \\
              0 &  0
          \end{array} \right]=(P|T^2|^2P)^{\frac{1}{2}}\geq
          P|T^2|P\geq P|T^{\ast}|^2P $$$$=PTT^{\ast}P=\left[
         \begin{array}{rr}
              M_{|\lambda|^2}+AA^{\ast} & 0 \\
              0 &  0
          \end{array} \right]
,$$ where $A=PM_{w}EM_{u}- M_{\lambda}$. This implies that $A=0$
and so $T^{\ast}(f)=M_{\bar{u}}EM_{\bar{w}}(f)=\bar{\lambda} f$.
This means
that $\sigma_{p}(M_wEM_u)=\sigma_{jp}(M_wEM_u)$.\\

\vspace*{0.3cm} {\bf Corollary 2.9.} If $T=M_wEM_u$ is
$p$-hyponormal and $S(E(u))=X$, then\\

(1) $\sigma_{p}(M_wEM_u)=\sigma_{jp}(M_wEM_u)$.\\

(2) $\sigma_{jp}(M_wEM_u)\setminus
\{0\}=\{\lambda\in\mathbb{C}:\mu(A_{\lambda,w})>0\}\setminus
\{0\}$,\\
where $A_{\lambda,w}=\{x\in X:E(uw)(x)=\lambda\}$.\\

\section{ \sc\bf Examples }

In this section we present some examples of conditional
expectations and corresponding multiplication operators to
illustrate concrete application of the main results of the paper.\\

 {\bf Example 3.1.} Let $X=[0,1]$, $\Sigma$=sigma algebra of Lebesgue
measurable subset of $X$, $\mu$=Lebesgue measure on $X$. Fix
$n\in\{2, 3, 4 . . .\}$ and let $s:[0, 1]\rightarrow[0, 1]$ be
defined by $s(x)=x+\frac{1}{n}$(mod 1). Let $\mathcal{B}=\{E\in
\Sigma:s^{-1}(E)=E\}$. In this case
$E^{\mathcal{B}}(f)(x)=\sum^{n-1}_{j=0}f(s^{j}(x))$, where $s^j$
denotes the jth iteration of $s$. The functions $f$ in the range
of $E^{\mathcal{B}}$ are those for which the $n$ graphs of $f$
restricted to the intervals $[\frac{j-1}{n},\frac{j}{n}]$, $1\leq
j\leq n$, are all congruent. Also, $|f|\leq nE^{\mathcal{B}}(|f|)$
a.e. Hence, the operator $EM_u$ is bounded on $L^2([0,1])$ if and
only if $u\in L^{\infty}([0,1])$. These
operators are closely related to averaging operators.\\

{\bf Example 3.2.} Let $X=[0,1]\times [0,1]$, $d\mu=dxdy$,
$\Sigma$ the Lebesgue subsets of $X$ and let
$\mathcal{A}=\{A\times [0,1]: A \ \mbox{is a Lebesgue set in} \
[0,1]\}$. Then, for each $f$ in $L^2(\Sigma)$, $(Ef)(x,
y)=\int_0^1f(x,t)dt$, which is independent of the second
coordinate. Now, if we take $u(x,y)=y^{\frac{x}{8}}$ and $w(x,
y)=\sqrt{(4+x)y}$, then $E(|u|^2)(x,y)=\frac{4}{4+x}$ and
$E(|w|^2)(x,y)=\frac{4+x}{2}$. So, $E(|u|^2)(x,y)E(|w|^2)(x,y)=2$
and $|E(uw)|^2(x,y)=64\frac{4+x}{(x+12)^2}$. Direct computations
shows that $E(|u|^2)(x,y)E(|w|^2)(x,y)\leq|E(uw)|^2(x,y)$, a.e,
and since the operator $M_wEM_u$ is bounded, then
$E(|u|^2)(x,y)E(|w|^2)(x,y)\geq|E(uw)|^2(x,y)$ a.e. Thus,
$E(|u|^2)(x,y)E(|w|^2)(x,y)=|E(uw)|^2(x,y)$, a.e. Therefore by
Theorem 2.8 we conclude that
$$\sigma_{p}(M_wEM_u)=\sigma_{jp}(M_wEM_u)=\{\sqrt{2}\}.$$\\

{\bf Example 3.3.} Let $X=[0,1)$, $\Sigma$=sigma algebra of
Lebesgue measurable subset of $X$, $\mu$=Lebesgue measure on $X$.
Let $s:[0, 1)\rightarrow[0, 1)$ be defined by
$s(x)=x+\frac{1}{2}$(mod 1). Let $\mathcal{B}=\{E\in
\Sigma:s(E)=E\}$. In this case
$$E^{\mathcal{B}}(f)(x)=\frac{f(x)+f(s(x))}{2}.$$  Also, $|f|\leq
E^{\mathcal{B}}(|f|)$ a.e. Hence, the operator $EM_u$ is bounded
on $L^2([0,1))$ if and only if $u\in L^{\infty}([0,1))$. Let
$u(x)=\sqrt{x}$ and $w\equiv1$, then $|E(uw)|^2\geq
E(|u|^2)E(|w|^2)$ a.e. Thus by Theorem 2.8 we have
$$\sigma_{p}(M_wEM_u)=\sigma_{jp}(M_wEM_u)=\emptyset.$$\\

{\bf Example 3.4.}  Let $X=\mathbb{N}$,
$\mathcal{G}=2^{\mathbb{N}}$ and let $\mu(\{x\})=pq^{x-1}$, for
each $x\in X$, $0\leq p\leq 1$ and $q=1-p$. Elementary
calculations show that $\mu$ is a probability measure on
$\mathcal{G}$. Let $\mathcal{A}$ be the $\sigma$-algebra generated
by the partition $B=\{X_1=\{3n:n\geq1\}, X^{c}_1\}$ of $X$. So,
for every $f\in \mathcal{D}(E^{\mathcal{A}})$,

$$E(f)=\alpha_1\chi_{X_1}+\alpha_2\chi_{X^c_1}$$
and direct computations show that

$$\alpha_1=\frac{\sum_{n\geq1}f(3n)pq^{3n-1}}{\sum_{n\geq1}pq^{3n-1}}$$
and
$$\alpha_2=\frac{\sum_{n\geq1}f(n)pq^{n-1}-\sum_{n\geq1}f(3n)pq^{3n-1}}{\sum_{n\geq1}pq^{n-1}-\sum_{n\geq1}pq^{3n-1}}.$$

For example, if we set $f(x)=x$, then $E(f)$ is a special function
as follows;

$$\alpha_1=\frac{3}{1-q^3}, \ \ \ \ \
\alpha_2=\frac{1+q^6-3q^4+4q^3-3q^2}{(1-q^2)(1-q^3)}.$$\\
So, if $u$ and $w$ are real functions on $X$ such that $M_wEM_u$
is boundede on $l^p$, then
$$\sigma(M_wEM_u)=\sigma_p(M_wEM_u)=\{\frac{\sum_{n\geq1}u(3n)w(2n)pq^{3n-1}}{\sum_{n\geq1}pq^{3n-1}},
\frac{\sum_{n\geq1}u(n)w(n)pq^{n-1}-\sum_{n\geq1}u(3n)w(3n)pq^{3n-1}}{\sum_{n\geq1}pq^{n-1}-\sum_{n\geq1}pq^{3n-1}}\}.$$\\

%

\end{document}